\documentclass[12pt]{amsart}
\usepackage[all]{xy}
\usepackage{amssymb,latexsym}

\numberwithin{equation}{section}

\newtheorem{theorem}{Theorem}[section]
\newtheorem{proposition}[theorem]{Proposition}

\newtheorem{corollary}[theorem]{Corollary}
\newtheorem{lemma}[theorem]{Lemma}

\theoremstyle{definition}

\newtheorem{remark}[theorem]{Remark}

\newtheorem{definition}[theorem]{Definition}

\let\oldmarginpar\marginpar
\renewcommand\marginpar[1]{\-\oldmarginpar[\raggedleft\small\sf
#1]{\raggedright\small\sf #1}}

\newcommand{\Z}{\mathbb{Z}}

\newcommand{\Q}{\mathbb{Q}}
\newcommand{\C}{\mathbb{C}}

\newcommand{\ca}{{\mathcal A}}

\newcommand{\cF}{{\mathcal F}}

\newcommand{\co}{{\mathcal O}}

\newcommand{\TT}{{\mathbb T}}
\newcommand{\dd}{{\mathbf d}}
\newcommand{\ee}{{\mathbf e}}

\newcommand{\sgn}{\operatorname{sgn}}

\newcommand{\rem}[1]{\left\langle#1\right\rangle}

\newcommand{\overunder}[2]{\!\begin{array}{c}\scriptstyle{#1}\\[-.1in]
-\!\!\!-\!\!\!-\\[-.1in]\scriptstyle{#2}
\end{array}\!}

\begin{document}

\title[Laurent expansions via quiver representations]
{Laurent expansions in cluster algebras via quiver representations}

\author{Philippe Caldero}
\address{Institut Camille Jordan\\ Universit\'e Claude Bernard Lyon I\\
69622 Villeurbanne Cedex, France}
\email{caldero@igd.univ-lyon1.fr}

\author{Andrei Zelevinsky}
\address{\noindent Department of Mathematics, Northeastern University,
Boston, MA 02115, USA}
\email{andrei@neu.edu}



\begin{abstract}
We study Laurent expansions of cluster variables in a cluster
algebra of rank~$2$ associated to a generalized Kronecker quiver.
In the case of the ordinary Kronecker quiver, we obtain explicit
expressions for Laurent expansions of the elements of the
canonical basis for the corresponding cluster algebra.
\end{abstract}

\date{April 6, 2006; revised May 1, 2006}

\dedicatory{To Alexander Alexandrovich Kirillov on the occasion of his seventieth birthday}

 \thanks{Andrei Zelevinsky's research
supported by NSF (DMS) grant \# 0500534
and by a Humboldt Research Award.}

\maketitle

\vspace{-.2in}


\section{Introduction}
Cluster algebras introduced in
\cite{fominzelevinsky1} have found
applications in a diverse variety of settings which include
(in no particular order) total positivity, Lie theory, quiver
representations, Teichm\"uller theory, Poisson geometry, discrete
dynamical systems, tropical geometry, and algebraic combinatorics.
See, e.g., \cite{caldkell, fg-survey, fomin-reading, FZ-conf}
and references therein.

Among these connections, the one with quiver representations has
been developed especially actively.
This development started with an observation made in \cite{marshreizel}
that the underlying combinatorial structure for a cluster algebra
has a natural interpretation in terms of quiver representations.
The subsequent work aimed to extend this interpretation from combinatorics to algebraic
properties of cluster algebras.
In the process, new concepts of cluster categories and
cluster-tilted algebras have been introduced and studied in
\cite{BMRRT,BMR} and many subsequent publications.
These new concepts extend the classical theory of quiver
representations and provide an interesting generalization of
classical tilting theory.

In this paper, we focus on one important algebraic feature
of cluster algebras: the \emph{Laurent phenomenon} established
in \cite{fominzelevinsky1}.
We will deal only with \emph{coefficient-free} cluster algebras.
In the nutshell, such an algebra is a commutative ring~$\ca$ (in fact, an integral domain) with a
family of distinguished generators (\emph{cluster variables}) grouped
into (overlapping) \emph{clusters} of the same finite
cardinality~$n$.
Each cluster is algebraically independent and generates the field
of fractions of~$\ca$.
Thus, every cluster variable can be uniquely expressed as a
rational function of the elements of every given cluster.
The Laurent phenomenon asserts that these rational functions are
in fact Laurent polynomials with integer coefficients.

We would like to know more about the coefficients of
these Laurent polynomials.
As conjectured by S.~Fomin and A.~Zelevinsky (see e.g. \cite{FZ-conf}),
these coefficients are positive integers; this is still proved only
in a few special cases.
We feel that a very promising tool in attacking this conjecture
is a geometric interpretation (found in \cite{caldchap} and generalized
in \cite{caldkell2}) of a coefficient in
question as the Euler-Poincar\'e characteristic of an appropriate
Grassmannian of quiver representations.
One of the main goals of this paper is to attract
attention to the problem (that we find very interesting) of studying
these Grassmannians and in
particular, finding an explicit way to compute their Euler-Poincar\'e
characteristics.

Our main new result is a complete solution of the latter problem
for the cluster algebra associated with a root system of affine
type $A_1^{(1)}$.
For this algebra, the positivity of the Laurent polynomials in question was
established by elementary means in \cite{sherzel}, while a
combinatorial expression for their coefficients was given in
\cite{musikerpropp}.
In this paper, we show that the Euler characteristic interpretation
implies an unbeatably simple explicit expression for every
coefficient in question as a product of two binomial coefficients.
After such an expression is found, proving it is not hard;
the point is however that this expression (left unnoticed in
\cite{sherzel,musikerpropp}) follows naturally from the geometric
study of the appropriate Grassmannians of quiver representations.
In this case, the underlying quiver is the Kronecker quiver
with two vertices and two arrows from one vertex to another.

We are happy to present these results in a paper dedicated to A.~A.~Kirillov.
One of the authors (A.Z.) had been very fortunate to have
A.~A.~Kirillov as one of his teachers at Moscow State University.
One of the most impressive features of Kirillov's teaching style
is his ability to explain mathematical ideas in the
simplest possible terms, clearing them of unnecessary technical
background so that they can be appreciated by inexperienced young researchers.
We are trying to follow his example in our exposition.
In particular, we never go beyond the ordinary quiver
representations; and we make the paper self-contained by
giving a new elementary proof of the
Euler characteristic interpretation for the generalized Kronecker quiver.

In Kirillov's spirit, we now state our new formula for
the cluster variables in type $A_1^{(1)}$ in a self-contained and
elementary way.
Let $x_1, x_2, x_3, \dots$ be a sequence of rational functions in
two independent variables $x_1, x_2$ defined recursively by
\begin{equation}
\label{eq:A11-recursion}
x_{n+1}=\frac{x_n^2+1}{x_{n-1}} \quad (n \geq 2).
\end{equation}
In Theorem~\ref{th:formula-xm-b=2} we show that, for every
$n \geq 0$, the term $x_{n+3}$ is given by
\begin{equation}
\label{eq:A11-formula-m-positive}
x_{n+3} = x_1^{-n-1} x_2^{-n} \left(x_2^{2(n+1)} +
\sum_{q + r \leq n} {n-r \choose q}{n+1-q \choose r} x_1^{2q}
x_2^{2r}\right);
\end{equation}
in particular, all terms of the sequence are Laurent polynomials
with positive coefficients in $x_1$ and $x_2$.

Note that if the exponent~$2$ in \eqref{eq:A11-recursion} gets
replaced by a positive integer $b > 2$, then each $x_n$ is still
an integer Laurent polynomial in $x_1$ and $x_2$ by the Laurent phenomenon
proved in \cite{fominzelevinsky1}.
However in this case no explicit combinatorial expression or
closed formula is known for the coefficients; even positivity of these coefficients is
still open.
Each of the coefficients in question is the Euler-Poincar\'e characteristic of an appropriate
Grassmannian of quiver representations for the generalized
Kronecker quiver~$Q_b$ with two vertices and~$b$ arrows from one vertex to another.
We find it a very interesting challenge to use this interpretation
for finding an explicit expression for the coefficients.

The paper is organized as follows.
In Section~\ref{sec:background} we provide some general background
on the Laurent phenomenon in cluster algebras and its
geometric interpretation.
In Section~\ref{sec:gen-kronecker} we provide a self-contained
treatment of the Laurent expansions of cluster variables in a cluster algebra of
rank~$2$ associated with the generalized Kronecker quiver~$Q_b$.
Our main technical tool that allows us to give a new proof for the
geometric interpretation of Laurent expansions
(Theorem~\ref{th:XM-cluster-Qb}) are the functors $T^+$ and $T^-$
on the  category of $Q_b$-representations (see Definition~\ref{def:T+-})
obtained by a slight modification of reflection functors in \cite{bgp}.
In Section~\ref{sec:kronecker} we work with the indecomposable
preprojective and preinjective representations of the classical
Kronecker quiver~$Q_2$ and prove equality \eqref{eq:A11-formula-m-positive} (Theorem~\ref{th:formula-xm-b=2}).
Finally, in Section~\ref{sec:kronecker-tube} we use the regular
indecomposable $Q_2$-representations to obtain
explicit Laurent expansions of the elements of the canonical basis
(constructed in \cite{sherzel}) in the corresponding cluster algebra.

\section{Some background}
\label{sec:background}

In this section we recall some background and results from
\cite{fominzelevinsky1,caldkell2}.
The definitions below are not the most general ones:
we will deal only with coefficient-free cluster algebras having
skew-symmetric exchange matrices (instead of more general
skew-symmetrizable ones).

Let $\cF = \Q(x_1, \dots, x_n)$ be the field of rational functions
in~$n$ independent variables.
Let~$B$ be a skew-symmetric integer $n \times n$ matrix.
We will associate to~$B$ a commutative subring $\ca(B) \subset \cF$ called the
(coefficient-free) \emph{cluster algebra}.
Let~$\TT_n$ be the \emph{$n$-regular tree} whose edges are labeled by the numbers $1, \dots, n$,
so that the $n$ edges emanating from each vertex receive different labels.
We write $t \overunder{k}{} t'$ to indicate that vertices
$t,t'\in\TT_n$ are joined by an edge labeled by~$k$.
We also fix some vertex $t_0 \in \TT_n$ and refer to $t_0$ as the
\emph{initial vertex}.
We associate to~$B$ and to every $t \in \TT_n$ a skew-symmetric integer
$n \times n$ matrix $B_t$, and an $n$-tuple $(x_{1;t}, \dots, x_{n;t})$
of elements of~$\cF$.
They are uniquely determined by the \emph{initial conditions}
$$B_{t_0} = B, \quad x_{j;t_0} = x_j,$$
and the \emph{mutation relations} given as follows.
Whenever $t \overunder{k}{} t'$, the matrices
$B_t = (b_{ij})$ and $B_{t'} = (b'_{ij})$ are related by
\begin{equation}
\label{eq:matrix-mutation}
b'_{ij} =
\begin{cases}
-b_{ij} & \text{if $i=k$ or $j=k$;} \\[.05in]
b_{ij} + \sgn(b_{ik}) \ [b_{ik}b_{kj}]_+
 & \text{otherwise,}
\end{cases}
\end{equation}
where we use the notation
\begin{align*}
[b]_+ &= \max(b,0); \\
\sgn(b) &=
\begin{cases}
-1 & \text{if $b<0$;}\\
0  & \text{if $b=0$;}\\
 1 & \text{if $b>0$;}
\end{cases}
\end{align*}
and we have $x_{j;t'} = x_{j;t}$ for $j \neq k$, while
$x_{k;t'}$ and $x_{k;t}$ satisfy the \emph{exchange relation}
\begin{equation}
\label{eq:exchange-rel-xx}
x_{k;t'}x_{k;t} = \ \prod_i x_{i;t}^{[b_{ik}]_+}
+ \ \prod_i x_{i;t}^{[-b_{ik}]_+} \, .
\end{equation}

We refer to $B_t$ as the \emph{exchange
matrix} at~$t$, and to $(x_{1;t}, \dots, x_{n;t})$ as
the \emph{cluster} at~$t$.
The elements $x_{j;t}$ are called \emph{cluster variables}
(note that there may be some equalities among them).
The cluster algebra $\ca(B)$ is defined as the subring of $\cF$
generated by all cluster variables.

Being an element of $\cF$, every cluster variable is a rational function
in $x_1, \dots, x_n$.
The following result sharpens this considerably.

\begin{theorem}[Laurent phenomenon {\cite[Theorem~3.1]{fominzelevinsky1}}]
\label{th:Laurent-phenom}
The cluster algebra~$\ca(B)$  is contained in the
Laurent polynomial ring $\Z[x_1^{\pm 1}, \dots, x_n^{\pm 1}]$;
equivalently, every cluster variable is an integer
Laurent polynomial in $x_1, \dots, x_n$.
\end{theorem}

We now turn to a geometric interpretation of the Laurent
polynomials in Theorem~\ref{th:Laurent-phenom}.
First of all, we represent a skew-symmetric integer $n \times n$ matrix $B = (b_{ij})$
by means of the quiver $Q = Q(B)$ with vertices
$[1,n] = \{1, \dots, n\}$, and $b_{ij}$ arrows from $i$ to $j$
whenever $b_{ij} > 0$ (thus, $Q$ is allowed to have multiple edges).
Following \cite{ca3}, we say that~$B$ is \emph{acyclic} if $Q(B)$
has no oriented cycles.
(Equivalently, $B$ is acyclic if and only if,
by a simultaneous permutation of rows and columns, we can make
$b_{ij} \geq 0$ for all $i > j$.)
In the rest of the section we assume that the initial exchange
matrix~$B$ is acyclic.

Let us recall some basics on quiver representations (we do not
attempt to give a self-contained introduction to the subject, just
fix some terminology and notation; for more information see e.g.,
\cite{dlabringel}).
Recall that a \emph{representation}~$M$ (over some field~$K$) of a quiver
$Q = Q(B)$ is given by assigning a finite-dimensional $K$-vector space
$M_i$ to every vertex~$i$ of~$Q$, and a $b$-tuple
$(\varphi_{ji}^{(1)}, \dots, \varphi_{ji}^{(b)})$ of
linear maps $M_i \to M_j$ to every arrow $i \to j$
of multiplicity~$b = b_{ij}$ in~$Q$.
For our current purposes it is sufficient to work over $K = \C$.
The \emph{morphisms} between quiver representations are defined in
a natural way, giving rise to the category of quiver
representations.
This category is abelian, hence there is a well defined
notion of \emph{indecomposable} representations.

The \emph{dimension} of a representation~$M$ is an integer vector
$\dd = (d_1, \dots, d_n)$ given by $d_i = \dim M_i$.
A quiver representation~$M$ of dimension~$\dd$ is called
\emph{rigid} if a generic representation of dimension~$\dd$ is
isomorphic to~$M$; equivalently, $M$ has no nontrivial
self-extensions.

A \emph{subrepresentation}~$N$ of a representation~$M$
is specified by a collection of subspaces $N_i \subset M_i$
such that $\varphi_{ji}^{(k)}(N_i) \subset N_j$
for all~$i$, $j$ and~$k$.
For a representation~$M$ of dimension~$\dd$, and any nonnegative integer vector
$\ee = (e_1, \dots, e_n)$ such that $e_i \leq d_i$ for all~$i$, let
${\rm Gr}_\ee(M)$ denote the variety of all subrepresentations
of dimension~$\ee$ in~$M$.
By the definition, ${\rm Gr}_\ee(M)$ is a closed subvariety in the
product of Grassmannians $\prod_i {\rm Gr}_{e_i}(M_i)$.

Let $\chi_\ee (M)$ denote the Euler-Poincar\'e characteristic
of ${\rm Gr}_\ee (M)$ (see e.g., \cite[Section~4.5]{fulton}).
We associate to any representation~$M$ of~$Q$
with dimension vector~$\dd$ a Laurent polynomial $X_M(x_1, \dots, x_n)$
given by
\begin{equation}
\label{eq:X-M}
X_M(x_1, \dots, x_n) = x_1^{- d_1} \cdots x_n^{- d_n}
\sum_\ee \chi_\ee (M) \prod_{i,j} (x_i^{d_j - e_j}
x_j^{e_i})^{[b_{ij}]_+};
\end{equation}
this is easily seen to be equivalent to the definition in \cite{caldchap}.

Now we are ready to state the following result obtained in
\cite[Theorem~3]{caldkell2}; it generalizes \cite[Theorem~3.4]{caldchap}.

\begin{theorem}
\label{th:XM-theorem}
Let~$\ca(B)$ be the (coefficient-free) cluster algebra with an acyclic
skew-symmetric initial exchange matrix~$B$.
The correspondence $M \mapsto X_M(x_1, \dots, x_n)$ is a bijection
between the set of isomorphism classes of indecomposable rigid
representations of the quiver~$Q(B)$, and the set of all cluster
variables in~$\ca(B)$ not belonging to the initial cluster $\{x_1, \dots, x_n\}$.
\end{theorem}

As shown in \cite[Corollary~1]{caldkell2}, for every
representation~$M$ of an acyclic quiver~$Q$, the dimension vector
$\dd =(d_1, \dots, d_n)$ of~$M$ is the \emph{denominator vector}
of the Laurent polynomial $X_M(x_1, \dots, x_n)$; that is, for every
$j \in [1,n]$, the minimum of exponents of $x_j$ in all the
monomials of $X_M$ is equal to $-d_j$.
The dimension vectors of indecomposable rigid representations are
called \emph{real Schur roots}.
This terminology comes from the well-known results due to V.~Kac \cite{kac}.
Namely, let $A=(a_{ij})$ be the \emph{Cartan counterpart} of~$B$,
that is, the symmetric integer $n \times n$ matrix with all diagonal
entries equal to~$2$, and off-diagonal entries given by
$a_{ij} = -|b_{ij}|$.
Then the dimension vectors of indecomposable representations of
$Q(B)$ are precisely the positive roots of the root system
associated to~$A$, expanded in the basis $\{\alpha_1, \dots, \alpha_n\}$
of simple roots.
Furthermore, a positive root $\alpha = \sum d_i \alpha_i$ is
\emph{real} if and only if there is a unique isomorphism class  of
indecomposable representations with dimension vector $(d_1, \dots, d_n)$.
We see in particular, that every real Schur root is a positive real root.
Note that the set of real Schur roots depends on the orientation of $Q(B)$,
in contrast with positive roots and with real positive roots.
There seems to be no easy way to distinguish real Schur roots
among all positive real roots, see \cite{derksenweyman}.

Returning to Theorem~\ref{th:XM-theorem}, we have the following corollary.

\begin{corollary}
\label{cor:denominator}
In the situation of Theorem~\ref{th:XM-theorem}, a cluster
variable in~$\ca(B)$ is uniquely determined by the denominator
vector in its Laurent expansion with respect to the initial
cluster.
Furthermore, the denominator vectors of the cluster
variables not belonging to the initial cluster are precisely the
real Schur roots of~$Q( B)$.
\end{corollary}

\section{Rank $2$ cluster algebras and generalized Kronecker quiver}
\label{sec:gen-kronecker}

In this section we discuss Theorem~\ref{th:XM-theorem} and
Corollary~\ref{cor:denominator} (and give their independent proofs)
for the cluster algebra~$\ca(B)$ associated with the matrix
$$B = \left(\begin{matrix} 0 & b\\
                           -b & 0
      \end{matrix}\right),$$
where~$b$ is a positive integer.
Unraveling the definitions in Section~\ref{sec:background}, we see
that $\ca(B)$ is the subring of the ambient field $\Q(x_1, x_2)$
generated by the cluster variables $x_m \,\, (m \in \Z)$ defined
recursively by the relations
\begin{equation}
\label{eq:rank2-b-recursion}
x_{m-1} x_{m+1} = x_m^b+1 \quad (m \in \Z)\ .
\end{equation}
The clusters are the pairs $\{x_m, x_{m+1}\}$ for all $m \in \Z$, and we
choose $\{x_1, x_2\}$ as the initial cluster.
Theorem~\ref{th:Laurent-phenom}
asserts that each $x_m$ is an integer Laurent polynomial in $x_1$ and $x_2$.
For $m \in \Z - \{1,2\}$, let $\alpha(m) \in \Z^2$ denote the
denominator vector of the cluster variable $x_m$.

If $b = 1$, an easy calculation (done in many places before) gives
$$x_3 = \frac{x_2 + 1}{x_1}, \,\, x_4 = \frac{x_1 + x_2 + 1}{x_1 x_2}, \,\,
x_5 = \frac{x_1 + 1}{x_2}, \,\, x_{m+5} = x_m \quad (m \in \Z)\ .$$
(As an easy exercise, one can check that these expressions agree
with Theorem~\ref{th:XM-theorem}.)
For the rest of this section we fix $b$ and assume that $b \geq 2$.

Let $S_{-1}(x), S_0(x), S_1(x), \dots$ be \emph{normalized Chebyshev polynomials of the second kind} defined
recursively by
\begin{equation}
\label{eq:chebind}
S_{-1}(x)=0, \quad S_{0}(x)=1, \quad
S_{n+1}(x)=x S_{n}(x)-S_{n-1}(x) \quad (n \geq 0).
\end{equation}
As an easy consequence of \eqref{eq:rank2-b-recursion}, the
denominator vectors $\alpha(m)$ of the cluster variables are given as follows.

\begin{proposition}
\label{pr:rank2-b-denominators}
For each $n \geq 0$, we have
\begin{equation}
\label{eq:denoms-Ab}
\alpha(n+3) = (S_{n}(b), S_{n-1}(b)), \quad
\alpha(-n) =  (S_{n-1}(b), S_{n}(b)) \ .
\end{equation}
\end{proposition}

We identify the lattice $\Z^2$ with the root lattice for the
Cartan matrix
$$A = \left(\begin{matrix} 2 & -b\\
                           -b & 2
      \end{matrix}\right),$$
by identifying the standard basis vectors $(1,0)$ and $(0,1)$ with the simple roots
$\alpha_1$ and $\alpha_2$ (for the properties of rank~$2$ root
systems see, e.g., \cite[Section~3.1]{sherzel}).
Under this identification, the denominator vectors in \eqref{eq:denoms-Ab}
are precisely the real positive roots.
This follows easily from the relations (for all $n \geq 0$)
\begin{equation}
\label{eq:s-alpha}
s_1 \alpha(-n) = \alpha(n+4) = \sigma \alpha(-n-1), \,\,
s_2 \alpha(n+3) = \alpha(-n-1) = \sigma \alpha(n+4),
\end{equation}
where $s_1$ and $s_2$ are simple reflections acting on $\Z^2$ by
matrices
$$s_1 = \left(\begin{matrix} -1 & b\\
                           0 & 1
      \end{matrix}\right),\quad
s_2 = \left(\begin{matrix} 1 & 0\\
                           b & -1
      \end{matrix}\right),$$
and $\sigma: \Z^2 \to \Z^2$ acts by interchanging the two
components of a vector (cf. \cite[(3.1)]{sherzel}).

Turning to Theorem~\ref{th:XM-theorem}, we note that the quiver
$Q_b = Q(B)$ associated to~$B$ has two vertices $1$ and $2$, and $b$
arrows from $1$ to $2$.
When $b=2$ (resp. $b > 2$), this quiver is called the
\emph{Kronecker quiver} (resp. \emph{generalized Kronecker quiver}).

A $Q_b$-representation $M$ over a field $K$ consists of a pair of finite dimensional $K$-vector spaces
$(M_1,M_2)$ and a $b$-tuple of linear maps $(\varphi_1,\ldots,\varphi_b)$ from $M_1$ to $M_2$.
We will use the following two interpretations of the tuple
$(\varphi_1,\ldots,\varphi_b)$: the \emph{column interpretation}
identifying it with a linear map $\varphi^c: M_1 \to M_2^b$, and
the \emph{row interpretation}
identifying it with a linear map $\varphi^r: M_1^b \to M_2$.

The dimension vector of $M$ is an integer vector ${\bf dim}\ M = (\dim M_1,\dim M_2) \in \Z^2$.
A subrepresentation $N$ of $M$ is a pair of subspaces $(N_1,N_2)$ such that
$N_i\subset M_i$ for $i=1,2$, and $\varphi_k(N_1)\subset N_2$ for $k \in [1,b]$.
In accordance with \eqref{eq:X-M}, we associate with every $Q_b$-representation
$M$ of dimension $\dd = (d_1, d_2)$ a Laurent polynomial
\begin{equation}
\label{eq:X-M-Qb}
X_M(x_1, x_2) = x_1^{- d_1} x_2^{- d_2}
\sum_{\ee = (e_1, e_2)} \chi_\ee (M) x_1^{b(d_2 - e_2)}
x_2^{be_1}.
\end{equation}

Since the vectors $\alpha(m)$ for $m \in \Z - \{1,2\}$ are real
positive roots, by Kac's theorem \cite{kac}, each of them is the
dimension vector of a unique (up to an isomorphism) indecomposable
$Q_b$-representation $M(m)$.
In particular, $M(3) = S_1$ is a simple representation of
dimension $\alpha_1$, and $M(0) = S_2$ is a simple representation of
dimension $\alpha_2$.
It is well known (see e.g., \cite{derksenweyman}) that each $\alpha(m)$ is a real Schur root, i.e.,
all $M(m)$ are rigid, but we will not use this fact.
The representations $M(-n)$ (resp. $M(n+3)$) for $n \geq 0$ are
also known as \emph{preprojective} (resp. \emph{preinjective})
indecomposable $Q_b$-representations.

\begin{theorem}
\label{th:XM-cluster-Qb}
For every $m \in \Z - \{1,2\}$, the cluster variable
$x_m$ is equal to $X_{M(m)}(x_1, x_2)$.
\end{theorem}

Our main tool in proving Theorem~\ref{th:XM-cluster-Qb} will be the following
functors on the category of $Q_b$-representations.

\begin{definition}
\label{def:T+-}
Let $M=(M_1,M_2;\varphi_1,\ldots,\varphi_b)$ be a $Q_b$-representation.
\begin{itemize}
\item
The duality functor $D$ sends $M$
to $M^* = (M_2^*, M_1^*; \varphi_1^*,\ldots,\varphi_b^*)$.
\item

The functor $T^+$ sends $M$ to
$M^+=( M_1^+, M_2^+;\varphi_1^+,\ldots,\varphi_b^+)$ given by
 $$M_1^+=M_2,\,\,  M_2^+=\hbox{coker}(\varphi^c: M_1\rightarrow
 M_2^b),$$
and ${(\varphi^+)}^r: (M_1^+){^b}=M_2^b\rightarrow  M_2^+$ being the natural projection.

\item

The functor $T^-$ sends $M$ to
$M^- = (M_1^-, M_2^-;\varphi_1^-,\ldots,\varphi_b^-)$ given by
 $$M_2^-=M_1,\,\,  M_1^-=\ker(\varphi^r: M_1^b\rightarrow M_2),$$
and ${(\varphi^-)}^c: M_1^-\rightarrow  M_1^b =(M_2^-){^b}$ being the natural embedding.
\end{itemize}
\end{definition}

All these functors are additive and send direct sums to direct sums.
The functors $T^+$ and $T^-$ are slight modifications of
\emph{reflection functors} from \cite{bgp}.
The following properties are immediate from the definition.

\begin{proposition}
\label{pr:compositions}
\begin{enumerate}
\item
$D^2 = {\rm Id}, \quad T^- = D T^+ D, \quad T^+ = D T^- D$.
\item
The composition $T^+ T^-$ sends a $Q_b$-representation
$M=(M_1,M_2;\varphi_1,\ldots,\varphi_b)$ to
$(M_1, M_1^b/\ker(\varphi^r); \psi_1,\ldots,\psi_b)$,
with $\psi^r$ being the natural projection $M_1^b \to
M_1^b/\ker(\varphi^r)$.
\item
The composition $T^- T^+$ sends a $Q_b$-representation
$M=(M_1,M_2;\varphi_1,\ldots,\varphi_b)$ to
$({\rm im}(\varphi^c), M_2; \psi_1,\ldots,\psi_b)$,
with $\psi^c$ being the natural embedding ${\rm im}(\varphi^c) \to M_2^b$.
\end{enumerate}
\end{proposition}

\begin{proposition}
\label{pr:T--image}
The following conditions on a $Q_b$-representation
$M$ are equivalent:
\begin{enumerate}
\item $M = T^- N$ for some representation $N$.
\item $M = T^- T^+ M$.
\item The map $\varphi^c: M_1 \to M_2^b$ is injective.
\item ${\bf dim}\ T^+ M = \sigma s_1 ({\bf dim}\ M)$.
\item $T^+ M' \neq 0$ for any non-zero subrepresentation $M'$ of $M$.
\item $S_1$ is not a direct summand of $M$.
\end{enumerate}
\end{proposition}

\begin{proof}
The implication $(2) \Longrightarrow (1)$ is trivial.
The implication $(1) \Longrightarrow (3)$ is immediate from the definition of $T^-$.
The equivalence $(2) \Longleftrightarrow (3)$ follows from
Proposition~\ref{pr:compositions} (3).
The equivalences $(3) \Longleftrightarrow (4) \Longleftrightarrow (5)$
are immediate from the definition of $T^+$.
The implication $(5) \Longrightarrow (6)$ is clear since $T^+ S_1 = 0$.
Finally, $(6) \Longrightarrow (3)$ is proved by contradiction as follows.
Suppose (3) does not hold, and choose a one-dimensional subspace
$M'_1 \subset \ker (\varphi^c)$.
Let $M''_1 \subset M_1$ be a subspace such that $M_1 = M'_1 \oplus M''_1$.
Then $M$ is a direct sum of subrepresentations $(M'_1,0)$ and
$(M''_1, M_2)$, and $(M'_1,0)$ is isomorphic to $S_1$, in
contradiction to (6).
\end{proof}

\begin{corollary}
\label{cor:T+indecomps}
Every indecomposable $Q_b$-representation $M$ not isomorphic to
$S_1$ satisfies equivalent conditions in
Proposition~\ref{pr:T--image}; furthermore, $T^+ M$ is also
indecomposable.
\end{corollary}

\begin{proof}
The first statement follows from condition (6) in
Proposition~\ref{pr:T--image}.
Now let $N = T^+M$, and suppose that $N$ is the direct sum of two
non-zero representations $N'$ and $N''$.
Applying the duality functor, we get $DN = DN' \oplus DN''$.
By Proposition~\ref{pr:compositions} (1), the representation
$DN = D T^+ M = T^- DM$ satisfies condition (1) in Proposition~\ref{pr:T--image}.
Therefore, by condition (5), both $T^+ DN'$ and $T^+ DN''$ are
non-zero.
Applying condition (2) and Proposition~\ref{pr:compositions} (1), we obtain
$M = T^- N = T^- N' \oplus T^- N'' =
D T^+ D N' \oplus D T^+ D N''$, in contradiction with the
assumption that $M$ is indecomposable.
\end{proof}

We now obtain an explicit description of
the indecomposable representations $M(m)$ for $m \in \Z - \{1,2\}$.

\begin{proposition}
\label{pr:Mm-description}
For every $n \geq 0$, we have
\begin{equation}
\label{eq:Mm-description}
M(-n) = (T^+){^n} S_2, \quad M(n+3) = (T^-){^n} S_1.
\end{equation}
\end{proposition}

\begin{proof}
Remembering the assumption $b \geq 2$, we start by observing that
all the roots $\alpha(m)$ are distinct, hence the corresponding indecomposable
representations $M(m)$ are mutually non-isomorphic.
In particular, all $M(-n)$ for $n \geq 0$ are not isomorphic to $S_1$.
To prove that $M(-n) = (T^+){^n} S_2$, we proceed by induction on~$n$.
The statement is trivial for $n=0$.
Now assume that it holds for some $n \geq 0$.
By Corollary~\ref{cor:T+indecomps}, the representation
$T^+ M(-n)$ is indecomposable.
Applying Proposition~\ref{pr:T--image} (4) and \eqref{eq:s-alpha},
we conclude that
$${\bf dim}\ T^+ M(-n) = \sigma s_1 ({\bf dim}\ M(-n)) =
\sigma s_1 (\alpha(-n)) = \alpha(-n-1).$$
Therefore, $T^+ M(-n) = M(-n-1)$, proving the first equality in
\eqref{eq:Mm-description}.
To prove the second equality in \eqref{eq:Mm-description} note that
$$M(n+3) = D M(-n) = D (T^+){^n} S_2 = (T^-){^n} D S_2 =
(T^-){^n} S_1.$$
\end{proof}

Turning to the proof of Theorem~\ref{th:XM-cluster-Qb}, we
start by rewriting \eqref{eq:X-M-Qb} as
\begin{equation}
\label{eq:XM-PM}
X_M(x_1, x_2) = x_1^{- d_1} x_2^{- d_2} P_M(x_1^b, x_2^b),
\end{equation}
where $P_M$ is a polynomial given by
\begin{equation}
\label{eq:PM-Qb}
P_M(z_1, z_2) =
\sum_{\ee = (e_1, e_2)} \chi_\ee (M) z_1^{d_2 - e_2} z_2^{e_1}.
\end{equation}

To work with $P_M(z_1, z_2)$, we need to recall some properties of the Euler-Poincar\'e characteristic.
We follow the treatment in \cite[Section~4.5]{fulton}, where
the Euler-Poincar\'e characteristic $\chi(X)$ is defined
for any complex algebraic variety~$X$ (not
necessarily smooth, projective or irreducible).
The following facts are shown in loc.cit.
\begin{align}
\label{eq:chi-affine}
&\text{If ${\mathbb A}$ is a finite dimensional affine space, then $\chi({\mathbb
A})=1$.}\\
\label{eq:chi-additivity}
&\text{If a variety $X$ is a disjoint union of finitely many}\\
\nonumber
&\text{locally closed subvarieties $X_i$,
then $\chi(X)= \sum \chi(X_i)$.}\\
\label{eq:chi-multiplicativity}
&\text{If $X\rightarrow  Z$ is a fiber bundle (locally trivial in the Zariski
topology)}\\
\nonumber
&\text{with fiber $Y$,
then $\chi(X)=\chi(Y)\chi(Z)$.}
\end{align}
As a consequence of \eqref{eq:chi-affine} and \eqref{eq:chi-additivity},
the Schubert cell decomposition of the Grassmannian implies that
\begin{equation}
\label{eq:chigras}
\chi({\rm Gr}_r(V))={\dim V\choose r}.
\end{equation}

Now let $M=(M_1,M_2;\varphi_1,\ldots,\varphi_b)$ be an arbitrary $Q_b$-representation
of dimension $(d_1,d_2)$.
For every two nonnegative integers $p$ and $r$, we set
\begin{align}
\label{eq:Zpr}
Z_{p,r}(M) &= \{U \in {\rm Gr}_r(M_1): \dim(\sum_{k=1}^b \varphi_k(U)) = d_2 - p\},\\
\label{eq:Z'pr}
Z'_{p,r}(M) &= \{U \in {\rm Gr}_{d_2-r}(M_2): \dim(\bigcap_{k=1}^b \varphi_k^{-1}(U)) = p\}.
\end{align}

\begin{proposition}
\label{pr:P-two-shifts}
We have
\begin{equation}
\label{eq:P-shifts}
P_M(z_1, z_2) = \sum_{p,r} \chi(Z_{p,r}(M)) (z_1+1)^p z_2^r
= \sum_{p,r} \chi(Z'_{p,r}(M)) z_1^r (z_2+1)^p.
\end{equation}
\end{proposition}

\begin{proof}
For every dimension vector
$\ee = (e_1, e_2)$, we split the Grassmannian ${\rm Gr}_{\ee}(M)$
into the disjoint union of subvarieties
\begin{equation}
\label{eq:GreM-p-1}
{\rm Gr}_{p,\ee}(M) = \{(N_1,N_2) \in {\rm Gr}_{\ee}(M):
N_1 \in Z_{p,e_1}(M)\}.
\end{equation}
The projection $(N_1, N_2) \mapsto N_1$ makes each
${\rm Gr}_{p,\ee}(M)$ into a fiber bundle over $Z_{p,e_1}(M)$.
Since, for a given $N_1 \in Z_{p,e_1}(M)$, the only condition on~$N_2$
is that $\sum_{k=1}^b \varphi_k(N_1) \subset N_2$,
the fiber of this bundle is the Grassmannian
of $(e_2 - d_2 + p)$-dimensional subspaces in a
$p$-dimensional vector space.
Applying \eqref{eq:chi-additivity},
\eqref{eq:chi-multiplicativity} and \eqref{eq:chigras}, we obtain
$$\chi_\ee(M) = \sum_p \chi({\rm Gr}_{p,\ee}(M)) =
\sum_p \binom{p}{e_2 - d_2 + p} \chi(Z_{p,e_1}(M)).$$
Substituting this expression into \eqref{eq:PM-Qb}, we obtain
\begin{equation*}
P_M(z_1, z_2) =
\sum_{p, e_1, e_2} \binom{p}{e_2 - d_2 + p} \chi(Z_{p,e_1}(M)) z_1^{d_2 - e_2} z_2^{e_1}
 = \sum_{p, r} \chi(Z_{p,r}(M)) (z_1+1)^{p} z_2^{r},
\end{equation*}
proving the first equality in \eqref{eq:P-shifts}.
The second equality can be proved in a similar way.
Alternatively, it is easy to see that the correspondence
$U \mapsto U^\perp \subset M_1^*$ is an isomorphism
between $Z_{p,r}(M)$ and $Z'_{p,r}(DM)$; this implies
the second equality in \eqref{eq:P-shifts} in view of
the (easily proved) observation
\begin{equation}
\label{eq:PM-D}
P_{DM}(z_1, z_2) = P_{M}(z_2, z_1).
\end{equation}
\end{proof}

The key ingredient of the proof of Theorem~\ref{th:XM-cluster-Qb}
is the following proposition.

\begin{proposition}
\label{pr:T+P}
Suppose $M$ is a $Q_b$-representation of dimension $(d_1, d_2)$
satisfying equivalent conditions in Proposition~\ref{pr:T--image}.
Then we have
\begin{equation}
\label{eq:T+P}
P_{T^+ M}(z_1, z_2) =
(z_1+1)^{- d_1} z_2^{d_2} P_M\left(\frac{(z_1+1)^b}{z_2}, z_1\right) \ .
\end{equation}
\end{proposition}

\begin{proof}
Consider the representation $T^+ M = M^+$ as defined
in Definition~\ref{def:T+-}.
By Proposition~\ref{pr:T--image} (4), we have
${\bf dim}\ T^+ M = (d_2, bd_2 - d_1)$.
Since $M_1^+ = M_2$, the statement of the following lemma makes sense.

\begin{lemma}
\label{lem:Zqr-M+}
Under the condition in Proposition~\ref{pr:T+P},
the variety $Z'_{p,r}(M)$ is equal to $Z_{p+br-d_1,d_2-r}(T^+M)$.
\end{lemma}

\begin{proof}
It suffices to show the following: if $U \subset M_2$
belongs to $Z'_{p,r}(M)$ then
$$\dim(\sum_{k=1}^b \varphi_k^+(U)) = b(d_2-r) - p.$$
By Definition~\ref{def:T+-}, we have
$$\sum_{k=1}^b \varphi_k^+(U) = (\varphi^+){^r}(U^b) = U^b/(U^b
\cap {\rm im}(\varphi^c)),$$
hence
$$\dim(\sum_{k=1}^b \varphi_k^+(U)) = b(d_2-r) - \dim(U^b \cap {\rm im}(\varphi^c)).$$
Remembering the definition of~$\varphi^c$, we see that
$$U^b \cap {\rm im}(\varphi^c) = \varphi^c(\bigcap_{k=1}^b \varphi_k^{-1}(U)).$$
Since $\varphi^c$ is injective by Proposition~\ref{pr:T--image} (3),
we conclude that
$$\dim(U^b \cap {\rm im}(\varphi^c)) = \dim(\bigcap_{k=1}^b
\varphi_k^{-1}(U))= p,$$
finishing the proof of Lemma~\ref{lem:Zqr-M+}.
\end{proof}

To prove \eqref{eq:T+P}, it suffices to combine Lemma~\ref{lem:Zqr-M+}
and formula \eqref{eq:P-shifts}:
\begin{align*}
P_{T^+ M}(z_1, z_2) &=
\sum_{p,r} \chi(Z'_{p,r}(M)) (z_1+1)^{p+br-d_1} z_2^{d_2-r}\\
&= (z_1+1)^{-d_1} z_2^{d_2} \sum_{p,r} \chi(Z'_{p,r}(M))
(z_1+1)^p \left(\frac{(z_1+1)^b}{z_2}\right)^{r}\\
&= (z_1+1)^{- d_1} z_2^{d_2} P_M\left(\frac{(z_1+1)^b}{z_2},z_1\right),
\end{align*}
as claimed.
\end{proof}

\begin{proof}[Proof of Theorem~\ref{th:XM-cluster-Qb}]
To show that $x_{-n} = X_{M(-n)}(x_1, x_2)$ for
$n \geq 0$, we proceed by induction on~$n$.
The check for $n=0$ is straightforward.
Thus we assume that $x_{-n} = X_{M(-n)}(x_1, x_2)$ for
some $n \geq 0$, and will show that
$x_{-n-1} = X_{M(-n-1)}(x_1, x_2)$.
To see this, we first note that $x_{-n-1} = X_{M(-n)}(x_0, x_1)$
by the obvious symmetry of
the exchange relations \eqref{eq:rank2-b-recursion}.
Now we apply Proposition~\ref{pr:T+P} to
$M = M(-n)$ and $T^+ M = M(-n-1)$
(see Corollary~\ref{cor:T+indecomps} and
Proposition~\ref{pr:Mm-description}).
Using \eqref{eq:XM-PM}, \eqref{eq:T+P} and \eqref{eq:rank2-b-recursion},
we obtain
\begin{align*}
X_{M(-n-1)}(x_1, x_2) &=
x_1^{- d_2} x_2^{- bd_2 + d_1} P_{M(-n-1)}(x_1^b, x_2^b)\\
&= x_1^{- d_2} x_2^{d_1} (x_1^b+1)^{- d_1}
P_{M(-n)} \left(\left(\frac{x_1^b+1}{x_2}\right)^b, x_1^b\right)\\
&= X_{M(-n)}\left(\frac{x_1^b+1}{x_2}, x_1\right)
= X_{M(-n)}(x_0, x_1) = x_{-n-1},
\end{align*}
as desired.

To show that $x_{n+3} = X_{M(n+3)}(x_1, x_2)$
for $n \geq 0$, we note that by \eqref{eq:PM-D} the equality
$M(n+3) = DM(-n)$ implies that
$X_{M(n+3)}(x_1, x_2)= X_{M(-n)}(x_2, x_1)$.
On the other hand, an obvious symmetry of the relations
\eqref{eq:rank2-b-recursion} implies that the automorphism of
$\Q(x_1, x_2)$ that interchanges $x_1$ and $x_2$, sends $x_{-n}$
to $x_{n+3}$; therefore, the Laurent expansion of $x_{n+3}$
in $x_1$ and $x_2$ is also obtained from that of $x_{-n}$
by interchanging $x_1$ and $x_2$.
This completes the proof of Theorem~\ref{th:XM-cluster-Qb}.
\end{proof}

\section{Cluster variables associated with the Kronecker quiver}
\label{sec:kronecker}

In this section we sharpen the results in Section~\ref{sec:gen-kronecker}
in the special case $b=2$.
Thus we work with the cluster algebra~$\ca(B)$
associated with the matrix
\begin{equation}
\label{eq:matrix-b=2}
B = \left(\begin{matrix} 0 & 2\\
                           -2 & 0
      \end{matrix}\right).
\end{equation}
In this case, we obtain the following explicit expression for
cluster variables.

\begin{theorem}
\label{th:formula-xm-b=2}
In the cluster algebra associated to the matrix $B$ in
\eqref{eq:matrix-b=2}, the cluster variables are given by
\begin{align}
\label{eq:A11-formula-m-negative}
x_{-n} &= x_1^{-n} x_2^{-n-1} \left(x_1^{2(n+1)} +
\sum_{q + r \leq n} {n+1-r \choose q}{n-q \choose r} x_1^{2q}
x_2^{2r}\right),\\
\label{eq:A11-formula-m-positive-2}
x_{n+3} &= x_1^{-n-1} x_2^{-n} \left(x_2^{2(n+1)} +
\sum_{q + r \leq n} {n-r \choose q}{n+1-q \choose r} x_1^{2q}
x_2^{2r}\right)
\end{align}
for all $n \geq 0$.
\end{theorem}

\begin{proof}
First of all, an obvious symmetry of the relations
\eqref{eq:rank2-b-recursion} implies that the map $x_m \mapsto x_{3-m}$
extends to an automorphism of our cluster algebra.
Applying this automorphism to both sides of \eqref{eq:A11-formula-m-negative}
yields \eqref{eq:A11-formula-m-positive-2}, so
it suffices to prove \eqref{eq:A11-formula-m-negative}.

We use Theorem~\ref{th:XM-cluster-Qb} to express $x_{-n}$ in terms
of the representations of the Kronecker
quiver $Q_2$ consisting of two vertices $1$ and $2$ and two
arrows from $1$ to $2$.
As a consequence of \eqref{eq:chebind}, we have $S_n(2) = n+1$
for all $n \geq -1$.
Thus, by Proposition~\ref{pr:rank2-b-denominators}, $x_{-n}$
has the denominator vector $(n,n+1)$, and so
the corresponding indecomposable $Q_2$ representation $M(-n)$
is of dimension $(n,n+1)$ as well.
An explicit form of this representation can be given as follows.

\begin{proposition}
\label{pr:M-n-explicit-Q2}
Let $M(-n) = (M_1, M_2; \varphi_1, \varphi_2)$
be the indecomposable $Q_2$ representation of dimension $(n,n+1)$.
Then there exist a basis $\{u_1, \dots, u_n\}$ in $M_1$ and a
basis $\{v_1, \dots, v_{n+1}\}$ in $M_2$ such that
$\varphi_1(u_k) = v_k$ and $\varphi_2(u_k) = v_{k+1}$ for $k \in [1,n]$.
\end{proposition}

This result is due to L.~Kronecker \cite{kronecker}; for a modern
treatment see \cite[Section~5.4]{FMV}.
A self-contained proof can be given by induction on~$n$ with the
help of Proposition~\ref{pr:Mm-description}.

The key ingredient for the proof of
\eqref{eq:A11-formula-m-negative} is the following result.

\begin{proposition}
\label{pr:chi-A11-preprojective}
Let $M(-n)$ be a $Q_2$-representation in
Proposition~\ref{pr:M-n-explicit-Q2}.
For every nonnegative integers $p, r$, we have (see
\eqref{eq:Zpr})
\begin{equation}
\label{eq:chi-A11-preprojective}
\chi(Z_{p,r}(M(-n))) = \binom{r-1}{n-p-r}  \binom{n+1-r}{p}
\end{equation}
(with the convention that the right hand side is equal to
$\delta_{p,n+1}$ for $r=0$).
\end{proposition}

\begin{proof}
The statement is trivial for $r=0$, so we assume that $r > 0$.
We use the Schubert cell decomposition of the Grassmannian ${\rm Gr}_r(M_1)$.
We label the Schubert cells by $r$-element subsets $J \subset [1,n]$.
The elements of the cell $\co(J)$ are parameterized by the arrays
of complex numbers
$$C = (c_{ij}) \,\, (j \in J, \ i \in [1,n] - J, \ i < j),$$
with the corresponding $U(C) \in \co(J)$ being
an $r$-dimensional subspace of $M_1$ with the basis
$$\{u_j + \sum_i c_{ij} u_i:  j \in J\}.$$
Breaking the subvariety $Z_{p,r}(M(-n)) \subset {\rm Gr}_r(M_1)$
into the disjoint union of its intersections with the Schubert
cells, and using \eqref{eq:chi-additivity}, we see that
$$\chi(Z_{p,r}(M(-n))) = \sum_J \chi(Z_{p,r}(M(-n)) \cap \co(J)).$$
Let $c(J)$ denote the number of connected components of a subset
$J \in [1,n]$ (by a connected component of~$J$ we mean a maximal interval
$[a,b] = \{a, a+1, \dots, b\}$ contained in~$J$).
The desired formula is a consequence of the following two results:
\begin{equation}
\label{eq:chi-Z-cap-O}
\chi(Z_{p,r}(M(-n)) \cap \co(J)) = \delta_{c(J), n+1-p-r};
\end{equation}
\begin{align}
\label{eq:schubert-cell-count}
&\text{the number of $r$-element subsets $J \subset [1,n]$ with
$c(J) = t$}\\
\nonumber
&\text{is equal to $\binom{r-1}{t-1}  \binom{n+1-r}{t}$.}
\end{align}

Since \eqref{eq:schubert-cell-count} is a purely combinatorial
statement, let us dispose of it first.
Let us write an $r$-element subset $J$ as the union of its
connected components:
$$J = [a_1, b_1] \cup \cdots \cup [a_t, b_t],$$
so we have
$$1 \leq a_1 < b_1 + 1 < a_2 < b_2 + 1 < \cdots < a_t < b_t + 1
\leq n+1,$$
and
$$(b_1+ 1 - a_1) + \cdots + (b_t +1 - a_t) = r.$$
Let $b_k + 1 - a_k = r_k$ for $k = 1, \dots, t$.
The number of $t$-tuples $(r_1, \dots, r_t)$ of positive integers
with sum $r$ is known to be $\binom{r-1}{t-1}$ (these tuples are
in a bijection with $(t-1)$-element subsets of $[1,r-1]$ via
$(r_1, \dots, r_t) \mapsto \{r_1, r_1 + r_2, \dots, r_1 + \cdots + r_{t-1}\}$).
And for every given such tuple, the number of corresponding
subsets $J$ is equal to $\binom{n+1-r}{t}$: they are in a
bijection with $t$-element subsets of $[1,n+1-r]$ via
$$[a_1, b_1] \cup \cdots \cup [a_t, b_t] \mapsto
\{a_1, a_2 - r_1, \dots, a_t - r_1 - \cdots - r_{t-1}\}.$$
This proves \eqref{eq:schubert-cell-count}.

Turning to the proof of \eqref{eq:chi-Z-cap-O}, we restate it as follows.
Fix an $r$-element subset $J \subset [1,n]$ and break the Schubert cell
$\co(J)$ into the disjoint union of the fibers of the function
$d: \co(J) \to \Z_{\geq 0}$ given by
$$d(U(C)) = \dim (\varphi_1 (U(C)) + \varphi_2 (U(C))).$$
We need to show that $\chi(d^{-1}(r+c(J)) = 1$, while all the
other fibers have Euler characteristic~$0$.
By Proposition~\ref{pr:M-n-explicit-Q2}, the subspace
$\varphi_1 (U(C)) + \varphi_2 (U(C)) \subset M_2$ is spanned by
the vectors $\{e_j(C), e_j^+(C): j \in J\}$, where we use the
notation
$$e_j (C) = v_j + \sum_i c_{ij} v_i  \quad
e_j^+ (C) = v_{j+1} + \sum_i c_{ij} v_{i+1}.$$
Denote $J^+ = \{j+1: j \in J\} \subset [2,n+1]$.
Note that $J - J^+$ is a set of representatives of the connected
components of $J$, so $|J - J^+| = c(J)$.
Consider the set of spanning vectors
$$E(C) = \{e_j^+(C) \, (j \in J), \,\, e_j(C) \, (j \in J - J^+)\}$$
of cardinality $r + c(J)$.
The vectors from $E(C)$ are linearly independent
since they have distinct leading terms in the expansion in the basis
$v_1, \dots, v_{n+1}$ of $M_2$.
Using these leading terms, it is easy to see that for each
remaining spanning vector $e_j(C)$ with $j \in J \cap J^+$, there
is a unique vector of the form
$$e'_j(C) = \sum_{i < j, \ i \notin J \cup J^+}
c'_{ij} v_i$$
obtained from $e_j(C)$ by adding a linear combination of the
vectors $e_{j'}(C), e_{j'}^+(C) \in E(C)$ with $j' < j$.
Furthermore, each coefficient $c'_{ij}$ is of the form
\begin{equation}
\label{eq:reduced-cij}
c'_{ij} = c_{ij} + P_{ij},
\end{equation}
where $P_{ij}$ is a polynomial in the variables $c_{i',j'}$ for $i' < j' < j$.
Clearly, replacing each $e_j(C)$ for $j \in J \cap J^+$ by
$e'_j(C)$ does not change the rank of the collection
$\{e_j(C), e_j^+(C): j \in J\}$, hence we have
\begin{equation}
\label{eq:dUC}
d(U(C)) = r + c(J)+ {\rm rk} (e'_j(C): j \in J \cap J^+);
\end{equation}
in particular, we see that $r+c(J)$ is the minimal value of the
function~$d$ on $\co(J)$.

In more geometric terms, the above statements can be rephrased as follows.
For each $j \in J \cap J^+$ let $V(j)$ denote the coordinate
subspace of $M_2$ spanned by the vectors $v_i$ for $i < j, \, i \notin J \cup J^+$.
The correspondence $U(C) \mapsto (e'_j(C): j \in J \cap J^+)$ defines a map
$\pi: \co(J) \to \prod_{j \in J \cap J^+} V(j)$.
As a consequence of \eqref{eq:reduced-cij}, $\pi$ is a fiber
bundle with fibers being finite dimensional affine spaces.
By \eqref{eq:chi-affine} and \eqref{eq:chi-multiplicativity}, we
have $\chi(X) = \chi(\pi(X))$ for every subvariety $X \subset \co(J)$.
In particular, in view of \eqref{eq:dUC}, $\chi(d^{-1}(r+c(J)+s))$
is equal to the Euler characteristic of the subvariety of
$\prod_{j \in J \cap J^+} V(j)$ consisting of all collections of
vectors having rank~$s$.
For $s=0$, the latter subvariety is just one point, implying
$\chi(d^{-1}(r+c(J)) = 1$.
And for $s > 0$, the subvariety of rank~$s$ collections in
$\prod_{j \in J \cap J^+} V(j)$ has Euler characteristic~$0$ by
\eqref{eq:chi-multiplicativity} since it has an obvious free
$\C^*$-action, and $\chi(\C^*) = 0$.
This completes the proof of
Proposition~\ref{pr:chi-A11-preprojective}.
\end{proof}

To finish the proof of \eqref{eq:A11-formula-m-negative},
we first use Theorem~\ref{th:XM-cluster-Qb} and \eqref{eq:XM-PM} to obtain
$$x_{-n} = X_{M(-n)}(x_1, x_2) = x_1^{- n} x_2^{- n-1} P_{M(-n)}(x_1^2, x_2^2).$$
Using \eqref{eq:chi-A11-preprojective} and \eqref{eq:P-shifts}, we get
\begin{align}
\label{eq:x-n-Laurent-shifted}
x_{-n} & = x_1^{-n} x_2^{-n-1} \, ((x_1^2+1)^{n+1}\\
\nonumber
& + \sum_{p \geq 0,\ r \geq 1}
\binom{r-1}{n-p-r}  \binom{n+1-r}{p}
(x_1^2+1)^p x_2^{2r}).
\end{align}

The desired formula \eqref{eq:A11-formula-m-negative}
follows from \eqref{eq:x-n-Laurent-shifted} by elementary
manipulations with binomial coefficients.
Expanding the powers of $(x_1^2+1)$, we obtain
$$x_{-n}x_1^{n} x_2^{n+1} = \sum_{q \geq 0} \binom{n+1}{q} x_1^{2q}
+ \sum_{q \geq 0,\ r \geq 1} a_{q,r} x_1^{2q} x_2^{2r},$$
where the coefficients $a_{q,r}$ are given by
$$a_{q,r} = \sum_{p} \binom{p}{q}
\binom{r-1}{n-p-r}  \binom{n+1-r}{p}.$$
Using an obvious identity
$$\binom{p}{q} \binom{n+1-r}{p} = \binom{n+1-r-q}{p-q} \binom{n+1-r}{q},$$
we can rewrite the last sum as
$$a_{q,r} = \binom{n+1-r}{q}
\sum_{p} \binom{r-1}{n-p-r}\binom{n+1-r-q}{p-q}.$$
Using the well-known Vandermonde identity
$$\sum_{k} \binom{a}{k}\binom{b}{c-k} = \binom{a+b}{c},$$
we conclude that
$$a_{q,r} = \binom{n+1-r}{q} \binom{n-q}{r},$$
implying \eqref{eq:A11-formula-m-negative}.
This completes the proof of Theorem~\ref{th:formula-xm-b=2}.
\end{proof}

\section{Regular representations of the Kronecker quiver}
\label{sec:kronecker-tube}

In this section we complement Theorem~\ref{th:formula-xm-b=2} by computing
the Laurent polynomials $X_M$ associated with the regular
indecomposable representations of the Kronecker quiver $Q_2$, i.e.,
those whose dimension vectors are imaginary positive roots.
In our case, the dimension vectors in question are $(n,n)$ for
$n \geq 1$, and for each $n$, the indecomposable
representations up to isomorphism are parameterized by $\C P^1$
(see \cite{kac} or \cite{FMV}).
It is easy to see that all the regular indecomposable representations~$M$
of the same dimension $(n,n)$ have the same Laurent polynomial $X_M$.
To compute it, we choose the following representative
$M^{\rm reg}(n)$ (cf. Proposition~\ref{pr:M-n-explicit-Q2}).

\begin{definition}
\label{def:M-regular}
Let $M^{\rm reg}(n)$ be a $Q_2$-representation of dimension $(n,n)$ defined as follows:
the space $M_1$ (resp. $M_2$) has a basis $\{u_1, \dots, u_n\}$
(resp. $\{v_1, \dots, v_n\}$) such that
$\varphi_1(u_k) = v_k$ and $\varphi_2(u_k) = v_{k+1}$ for $k \in [1,n]$,
with the convention that $v_{n+1} = 0$.
We denote
\begin{equation}
\label{eq:sn}
s_n = X_{M^{\rm reg}(n)}(x_1, x_2).
\end{equation}
\end{definition}

We prove the following analogue of Theorem~\ref{th:formula-xm-b=2}.

\begin{theorem}
\label{th:formula-sn}
The Laurent polynomials $s_n$ are given by
\begin{equation}
\label{eq:sn-formula}
s_{n} = x_1^{-n} x_2^{-n} \sum_{q+r \leq n}
{n-r \choose q}{n-q \choose r} x_1^{2q} x_2^{2r}
\end{equation}
for all $n \geq 1$.
\end{theorem}

\begin{proof}
The proof follows that of Theorem~\ref{th:formula-xm-b=2}.
The analogue of Proposition~\ref{pr:chi-A11-preprojective} is as follows.

\begin{proposition}
\label{pr:chi-regular}
For every nonnegative integers $p, r$, we have (see
\eqref{eq:Zpr})
\begin{equation}
\label{eq:chi-regular}
\chi(Z_{p,r}(M^{\rm reg}(n))) = \binom{r}{n-p-r}\binom{n-r}{p}.
\end{equation}
\end{proposition}

\begin{proof}
The proof of Proposition~\ref{pr:chi-regular} follows that of
Proposition~\ref{pr:chi-A11-preprojective} almost verbatim with
obvious modifications coming from the fact that $\varphi_2(u_n) = 0$.
First, \eqref{eq:chi-Z-cap-O} gets replaced by
\begin{equation}
\label{eq:chi-Z-cap-O-regular}
\chi(Z_{p,r}(M^{\rm reg}(n)) \cap \co(J)) = \delta_{c(J) - \varepsilon(J), n-p-r},
\end{equation}
where we set
\begin{equation}
\label{eq:epsilon}
\varepsilon(J) =
\begin{cases}
1 & \text{if $n \in J$;}\\
0  & \text{if $n \notin J$.}
\end{cases}
\end{equation}

We then show the following analogue of \eqref{eq:schubert-cell-count}:
\begin{align}
\label{eq:schubert-cell-count-regular}
&\text{the number of $r$-element subsets $J \subset [1,n]$ with
$c(J) - \varepsilon(J) = t$}\\
\nonumber
&\text{is equal to $\binom{r}{t}  \binom{n-r}{t}$.}
\end{align}
This can be proved by a slight modification of the proof of
\eqref{eq:schubert-cell-count}.
Alternatively, one can deduce \eqref{eq:schubert-cell-count-regular}
from \eqref{eq:schubert-cell-count}  by the following simple argument.
Let $c(r,n,t)$ denote the number of subsets $J$ in
\eqref{eq:schubert-cell-count}, that is, the number of $r$-element subsets $J \subset [1,n]$ with
$c(J) = t$.
Then it is easy to see that the number of subsets~$J$ in \eqref{eq:schubert-cell-count-regular}
is equal to $c(r,n-1,t) + c(r,n,t+1) - c(r,n-1,t+1)$.
Using \eqref{eq:schubert-cell-count}, we see that the number in
question is equal to
$$\binom{r-1}{t-1}  \binom{n-r}{t} + \binom{r-1}{t}  \binom{n+1-r}{t+1}
- \binom{r-1}{t}  \binom{n-r}{t+1}$$
$$= \binom{r-1}{t-1} \binom{n-r}{t}
+ \binom{r-1}{t}  \binom{n-r}{t} = \binom{r}{t}  \binom{n-r}{t},$$
as desired.

Formula \eqref{eq:chi-regular} is an immediate consequence of
\eqref{eq:chi-Z-cap-O-regular} and
\eqref{eq:schubert-cell-count-regular}.
\end{proof}

Arguing as in Section~\ref{sec:kronecker}, we obtain the following
analogue of \eqref{eq:x-n-Laurent-shifted}:
\begin{equation}
\label{eq:s-n-Laurent-shifted}
s_{n} = \frac{1}{x_1^{n} x_2^{n}} \sum_{p, r \geq 0}
\binom{r}{n-p-r}  \binom{n-r}{p}
(x_1^2+1)^p x_2^{2r}.
\end{equation}
Formula \eqref{eq:sn-formula} follows from \eqref{eq:s-n-Laurent-shifted}
in the same way as \eqref{eq:A11-formula-m-negative}
follows from \eqref{eq:x-n-Laurent-shifted}.
\end{proof}

Theorem~\ref{th:formula-sn}
allows us to sharpen the results in \cite{sherzel} on the canonical basis in the
cluster algebra $\ca$ associated to the Kronecker quiver.
Following \cite{sherzel}, we call a non-zero element $x \in \ca$
\emph{positive} if its Laurent expansion in terms of every cluster
$\{x_m, x_{m+1}\}$ has positive (integer) coefficients;
furthermore, a positive element $x$ is \emph{indecomposable} if it
cannot be written as a sum of two positive elements.
In \cite[Theorem~2.3]{sherzel} it is proved that all
indecomposable positive elements form a $\Z$-basis in~$\ca$,
referred to as the \emph{canonical basis}.
As shown in \cite[Theorem~2.8]{sherzel}, the canonical basis
consists of all \emph{cluster monomials}
$x_m^p x_{m+1}^q \,\, (m \in \Z, \, p, q \geq 0)$ together
with a sequence of elements $z_n \,\, (n \geq 1)$ defined as follows.
First of all, let
\begin{equation}
\label{eq:z1}
z_1 = \frac{x_1^2 + x_2^2 + 1}{x_1 x_2};
\end{equation}
the fact that $z_1 \in \ca$ follows from an easily checked equality
$z_1 = x_0 x_3 - x_1 x_2$.
The elements $z_n$ for all $n \geq 1$ are defined by
\begin{equation}
\label{eq:zn}
z_n = P_n(z_1),
\end{equation}
where the $P_n$ are \emph{normalized Chebyshev polynomials of the first kind},
related to the polynomials $S_n(x)$ in \eqref{eq:chebind} by
\begin{equation}
\label{eq:Tn-Un}
P_{n}(x)= S_{n}(x)-S_{n-2}(x) \quad (n \geq 0)
\end{equation}
(with the convention that $S_{-2}(x) = 0$).

We can now state an explicit formula (unnoticed in \cite{sherzel})
for the Laurent expansion of each $z_n$.

\begin{theorem}
\label{th:formula-zn}
The Laurent expansion of each $z_n$ for $n \geq 1$ in terms of
$x_1$ and $x_2$ is given by
\begin{equation}
\label{eq:zn-formula}
z_{n}
= x_1^{-n} x_2^{-n} \, (x_1^{2n} + x_2^{2n}
+ \sum_{q+r \leq n-1} \frac{n}{n-q-r}
{n-1-r \choose q}{n-1-q \choose r} x_1^{2q} x_2^{2r}).
\end{equation}
\end{theorem}

\begin{proof}
First of all, in view of \eqref{eq:z1}, formula
\eqref{eq:zn-formula} holds for $n=1$, and we also have $z_1 = s_1$.
A direct calculation using \eqref{eq:sn-formula} shows that the right hand side of
\eqref{eq:zn-formula} is equal to $s_n - s_{n-2}$ for $n \geq 2$
(with the convention that $s_0 = 1$).
Taking into account \eqref{eq:zn} and \eqref{eq:Tn-Un}, we see
that it remains to show that
\begin{equation}
\label{eq:sn-Un}
s_n = S_n(z_1) \quad (n \geq 0).
\end{equation}
By \eqref{eq:chebind}, it suffices to show that the elements $s_n$
satisfy the recursion
$$s_{n+1} = z_1 s_n - s_{n-1} \quad (n \geq 1).$$
This is an easy consequence of \eqref{eq:sn-formula}, finishing
the proof.
\end{proof}

We conclude with three remarks.

\begin{remark}
Explicit expressions \eqref{eq:A11-formula-m-negative},
\eqref{eq:A11-formula-m-positive-2} and \eqref{eq:zn-formula}
make obvious the facts about the Newton polygons of the elements
of the canonical basis in \cite[Propositions~3.5, 5.1 and 5.2]{sherzel}.
\end{remark}

\begin{remark}
In view of \eqref{eq:zn}, \eqref{eq:Tn-Un} and \eqref{eq:sn-Un},
the elements $z_n$ and $s_n$ of the cluster algebra $\ca$ are related by
\begin{equation}
\label{eq:zn-sn}
z_1 = s_1, \,\, z_n = s_n - s_{n-2} \quad (n \geq 2).
\end{equation}
It follows that replacing each $z_n$ by $s_n$ transforms the canonical
basis into another $\Z$-basis of $\ca$.
The relationship between this new basis and the canonical basis is
analogous to the relationship between the (dual) semicanonical and
the (dual) canonical basis for quantum groups, cf. \cite{GLS}.
\end{remark}

\begin{remark}
In view of formula \eqref{eq:X-M-Qb}, Theorems~\ref{th:formula-xm-b=2} and
\ref{th:formula-sn} provide a simple closed expression for the
Euler-Poincar\'e characteristic $\chi({\rm Gr}_\ee(M))$
of each ``quiver Grassmannian" in an arbitrary indecomposable representation
$M$ of the Kronecker quiver.
One can also use the proofs of these theorems to obtain a nice
combinatorial interpretation for $\chi({\rm Gr}_\ee(M))$.
Namely, if we realize $M$ as in Proposition~\ref{pr:M-n-explicit-Q2}
and Definition~\ref{def:M-regular}, then in each case, the spaces $M_1$ and $M_2$
are supplied with the distinguished bases $\{u_i\}$ and $\{v_i\}$, respectively.
The calculations in the proofs of Theorems~\ref{th:formula-xm-b=2} and
\ref{th:formula-sn} imply that $\chi({\rm Gr}_\ee(M))$ is equal
to the (finite) number of points $(N_1, N_2) \in {\rm Gr}_\ee(M)$
such that  both $N_1$ and $N_2$ are coordinate
subspaces with respect to these distinguished bases.
This expression for $\chi({\rm Gr}_\ee(M))$ can be rephrased as the following
combinatorial expression for the polynomial $P_M(z_1, z_2)$ (see \eqref{eq:PM-Qb}).
Consider the polynomials $F(w_1, \dots, w_N)$ given by
\begin{equation}
\label{eq:Fw}
F(w_1, \dots, w_N) = \sum_D \prod_{k \in D} w_k,
\end{equation}
where $D$ runs over all subsets of $[1,N]$ containing no two consecutive integers
(these polynomials appear in a different context in
\cite[Example~2.15]{yga}).
Then, for every $n \geq 0$, we have
\begin{align}
\nonumber
&P_{M(-n)}(z_1, z_2) = F(w_1, \dots, w_{2n+1})|_{w_k = z_{\rem{k}}},\\
\label{eq:P-Fw}
&P_{M(n+3)}(z_1, z_2) = F(w_1, \dots, w_{2n+1})|_{w_k = z_{\rem{k+1}}},\\
\nonumber
&P_{M^{\rm reg}(n+1)}(z_1, z_2) = F(w_1, \dots, w_{2n+2})|_{w_k = z_{\rem{k}}},
\end{align}
where $\rem{k}$ stands for the element of $\{1,2\}$ congruent to $k$ modulo~$2$.
In view of \eqref{eq:XM-PM}, we also have
\begin{align}
\nonumber
&x_{-n} = x_1^{-n} x_2^{-n-1} F(w_1, \dots, w_{2n+1})|_{w_k = x_{\rem{k}}^2},\\
\label{eq:xs-Fw}
&x_{n+3} = x_1^{-n-1} x_2^{-n} F(w_1, \dots, w_{2n+1})|_{w_k = x_{\rem{k+1}}^2},\\
\nonumber
&s_{n+1} = x_1^{-n-1} x_2^{-n-1} F(w_1, \dots, w_{2n+2})|_{w_k = x_{\rem{k}}^2}.
\end{align}
These formulas are easily seen to be equivalent to the combinatorial expressions
for cluster variables in \cite{musikerpropp}.
\end{remark}

\section*{Acknowledgments}
Andrei Zelevinsky thanks Claus Michael Ringel and Andrew Hubery
for stimulating discussions and their hospitality during his stay
at the University of Bielefeld.


\begin{thebibliography}{99}



\bibitem{ca3}
A.~Berenstein, S.~Fomin, A.~Zelevinsky,
Cluster algebras~III: Upper bounds and double Bruhat cells,
\textsl{Duke Math.~J.} \textbf{126} (2005), 1--52.

\bibitem{bgp}
J. Bernstein, I. Gelfand, V. Ponomarev,
Coxeter functors, and Gabriel's theorem,
\textsl{Uspehi Mat. Nauk} \textbf{28} (1973), no. 2(170), 19–-33.


\bibitem{BMR}
A.~Buan, R.~Marsh, I. Reiten,
Cluster tilted algebras, \texttt{math.RT/0402075}.





\bibitem{BMRRT}
A.~Buan, R.~Marsh, M. Reineke, I. Reiten, G. Todorov,
Tilting theory and cluster combinatorics, \texttt{math.RT/0402054}.

\bibitem{caldchap}
P.~Caldero, F. Chapoton,
Cluster algebras as Hall algebras of quiver representations,
\texttt{math.RT/0410184}.



\bibitem{caldkell}
P.~Caldero, B. Keller,
From triangulated categories to cluster algebras,
\texttt{math.RT/0506018}.

\bibitem{caldkell2}
P.~Caldero, B. Keller,
From triangulated categories to cluster algebras II,
\texttt{math.RT/0510251}.

\bibitem{derksenweyman}
H.~Derksen, J.~Weyman,
On canonical decomposition for quiver representations,
\textsl{Compositio Math.} \textbf{133} (2002), no. 3, 245–-265.

\bibitem{dlabringel}
V.~Dlab, C.~M.~Ringel,
Indecomposable representations of graphs and algebras,
\textsl{Mem.  Amer. Math. Soc.} \textbf{6} (1976), no. 173.


\bibitem{fg-survey}
V.~Fock, A.~Goncharov,
Dual Teichmuller and lamination spaces,
\texttt{math.DG/0510312}.

\bibitem{fomin-reading}
S.~Fomin, N.~Reading,
Root systems and generalized associahedra,
Lecture notes for the IAS/Park City Graduate Summer School in
Geometric Combinatorics (July 2004), \texttt{math.CO/0505518}.


\bibitem{fominzelevinsky1}
 S. Fomin, A. Zelevinsky,
Cluster algebras. I. Foundations,
\textsl{J. Amer. Math. Soc.}  \textbf{15}  (2002),  no. 2, 497--529.

\bibitem{fominzelevinsky2}
S. Fomin, A. Zelevinsky,
Cluster algebras. II. Finite type classification,
\textsl{Invent. Math.}  \textbf{154}  (2003),  no. 1, 63--121.

\bibitem{yga}
 S.~Fomin, A.~Zelevinsky,
$Y$-systems and generalized associahedra,
\textsl{Ann.\ in Math.} \textbf{158} (2003), 977--1018.


\bibitem{FZ-conf}
 S. Fomin, A. Zelevinsky,
Cluster algebras: Notes for the CDM-03 conference,
in: \textsl{CDM 2003: Current Developments in Mathematics,}
International Press, 2004.


\bibitem{FMV}
I. Frenkel, A. Malkin, M. Vybornov,
Affine Lie algebras and tame quivers,
\textsl{Selecta Math. (N.S.)}  \textbf{7}  (2001),  no. 1, 1--56.



\bibitem{fulton}
W. Fulton, Introduction to Toric Varieties,
Ann. of Math. Stud., vol. 131, Princeton University Press, 1993.

\bibitem{GLS}
C. Geiss, B. Leclerc, J. Schr\"oer,
Semicanonical bases and preprojective algebras,
\textsl{Ann. Sci. \'Ecole Norm. Sup. (4)}  \textbf{38}  (2005),  no. 2, 193--253.

\bibitem{kac}
V.~Kac,
Infinite root systems, representations of graphs and invariant theory,
\textsl{Invent. Math.} \textbf{56} (1980), no. 1, 57–-92.


\bibitem{kronecker}
L. Kronecker,
Algebraische reduction der schaaren bilinearer formen,
\textsl{Sitzungsberichte Akad. Berlin} (1890), 1225–-1237.



\bibitem{marshreizel}
R.~Marsh, M.~Reineke, A.~Zelevinsky,
Generalized associahedra via quiver representations.
\textsl{Trans. Amer. Math. Soc.}  \textbf{355} (2003),  no. 10, 4171--4186.

\bibitem{musikerpropp}
G.~Musiker, J.~Propp,
Combinatorial interpretations for rank-two cluster algebras of affine type,
\texttt{math.CO/0602408}.


\bibitem{sherzel}
P. Sherman, A. Zelevinsky,
Positivity and canonical bases in rank 2 cluster algebras of finite and affine types,
\textsl{Moscow Math. J.}  \textbf{4}  (2004),  no. 4, 947--974.

\end{thebibliography}
\end{document}